\documentclass[12pt,a4paper]{article}

\usepackage{comment,xcolor}
\usepackage{amssymb,amsmath,amsthm}
\usepackage{graphicx}
\usepackage[shortlabels]{enumitem}

\newtheorem{theorem}{Theorem}

\newtheorem{lemma}[theorem]{Lemma}

\newcommand{\ol}[1]{\overline{#1}}

\begin{document}
\normalsize

\title{{\bf \Huge Lines in the plane with the $L_1$ metric}}

\pagestyle{myheadings} \markright{{\small{\sc I.~Kantor:  Lines in the plane with the $L_1$ metric  }}}  

\author{Ida Kantor%
		\thanks{Supported by grant 19-04113 of the Czech Science Foundation (GACR) and by Charles University project UNCE/SCI/004.}}

\date{
 Charles University, Prague\\
(e-mail: \texttt{ida@iuuk.mff.cuni.cz})
}


\maketitle

\renewcommand{\thefootnote}{\empty}

\maketitle

\begin{abstract}
A well-known theorem in plane geometry states that any set of $n$ non-collinear points in the plane determines at least $n$ lines. Chen and Chv\'{a}tal asked whether an analogous statement holds within the framework of finite metric spaces, with lines defined using the notion of {\em betweenness}. 

In this paper, we prove that in the plane with the $L_1$ (also called Manhattan) metric, a non-collinear set of $n$ points induces at least $\lceil n/2\rceil$ lines. This is an improvement of the previous lower bound of $n/37$, with substantially different proof. As a consequence, we also get the same lower bound for non-collinear point sets in the plane with the $L_{\infty}$ metric. 
\end{abstract}

%%%%%%%%%%%%%%%%%%%%%%%%%%%%%%%5
\section{Lines in finite metric spaces}

A well-known theorem in plane geometry states that $n$ points in the plane are either collinear, or they induce at least $n$ lines. Erd\H{o}s noticed in~\cite{erdHos1944problem} that this is a corollary of the Sylvester-Gallai theorem (which states that for any non-collinear set $X$ of points in the plane, some line passes through exactly two points of $X$). Also, it is a special case of a theorem proved by de Bruijn and Erd\H{o}s~(\cite{deBE48}) about incidence structures. 

The proof of de Bruijn and Erd\H{o}s involves neither measurements of distances, nor measurements of angles. As such it is part of {\em ordered geometry}, which revolves around the ternary notion of {\em betweenness}: a point $z$ is between $x$ and $y$ if it is an interior point of the line segment with endpoints $x$ and $y$. We will write $[xzy]$ to indicate that $z$ is between $x$ and $y$. In terms of Euclidean distance $\rho$, 

\vspace{2mm}
\begin{center}
$z$ is between $x$ and $y$ $\ $ $\Longleftrightarrow$ $\ $ $x,y,z$ are distinct points and $\rho(x,y)=\rho(x,z)+\rho(z,y)$. 
\end{center}
\vspace{2mm}

In an arbitrary metric space, this notion becomes {\em metric betweenness}, introduced in~\cite{Menger28}.  The concept of a line is also generalized quite naturally:

\vspace{2mm}
\begin{center}
{\em A {\em line $\ol{uv}$} consists of $u$ and $v$ and all points $w$ such that one of $u,v,w$ is between the other two.}
\end{center}
\vspace{2mm}

We say that the pair $\{u,v\}$ {\em induces} the line $\ol{uv}$. A line containing all of $X$ is called a {\em universal line}. If $X$ has an universal line, we call $X$ a {\em collinear set}.

Having defined lines in arbitrary metric spaces in 2006 (\cite{ChenChvatal08}), Chen and Chvátal asked which properties of lines in Euclidean space translate into this setting. In particular, they posed the following question:

\vspace{5mm}
\begin{centering}
{\em Is it true that every finite metric space $(X,\rho)$ induces at least $|X|$ lines, or there is a universal line?}
\end{centering}
\vspace{5mm} 

The Chen--Chv\'{a}tal question is still open, although a number of interesting results 
related to it have been proved; these are surveyed in~\cite{ChvSurvey}. Among them, 
let us mention for future reference a theorem of Aboulker, Chen, 
Huzhang, Kapadia, and Supko~(\cite{ACHKS}, Theorem 3.1): In an arbitrary metric 
space, every non-collinear set of $n$ points induces $\Omega(\sqrt{n})$ lines.

\begin{figure}
  \centering
    \includegraphics[width=0.6\textwidth]{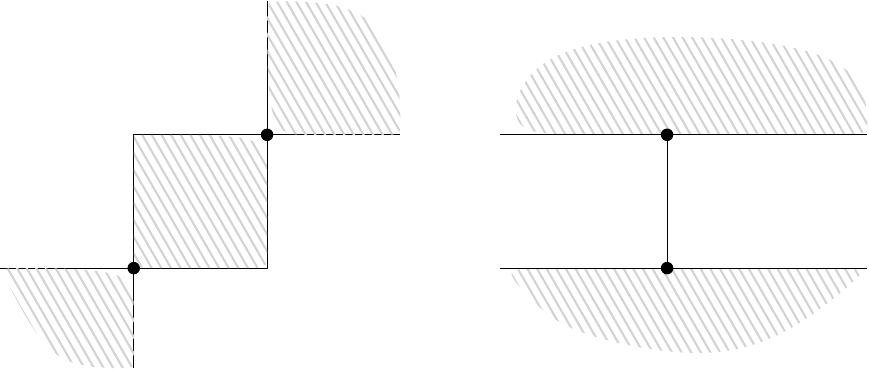}
		\caption{Line induced by an increasing pair and by a vertical pair.}
		\label{fig:twolines}
\end{figure}

In this paper, we concentrate on the plane with the $L_1$ (also called 
Manhattan) metric, defined by $\rho((x_1,y_1),(x_2,y_2))=|x_1-x_2|+|y_1-y_2|$. We encounter two very different kinds of lines in this case: those induced by two points that share a coordinate (horizontal or vertical pairs), and those induced by pairs that do not share a coordinate (so called increasing or decreasing pairs). The two types of lines (for increasing/decreasing pairs of points, and for horizontal/vertical pairs) can be seen in Figure~\ref{fig:twolines}. In this picture, the line in question is the finite subset of $X$ consisting of points that belong in the shaded area.

Bal\'{a}zs Patk\'{o}s and this author proved in~\cite{KP13} that every 
non-collinear set of $n$ points in this particular metric space induces 
at least $n/37$ lines. Moreover, if no two of the points share their $x$- 
or $y$-coordinate, then there are at least $n$ lines, i.e., the Chen---Chv\'{a}tal question has an affirmative answer. The proof relies on a lemma 
that was later extended and used by Aboulker, Chen, Huzhang, Kapadia, 
and Supko in deriving the weaker lower bound $\Omega(\sqrt{n})$ valid for all 
metric spaces.

In the present paper, we improve the lower bound of $n/37$ by a 
completely different method.

\begin{theorem}\label{main}
Let $X$ be a set of $n$ points in the plane with the $L_1$ metric. If there is no universal line, then $X$ induces at least $\lceil n/2\rceil$ lines.
\end{theorem}

Let us briefly mention the main ideas of the proof. Introduce an oriented graph $G$ on the vertex set $X$ as follows. Fix a line $L$ (this is a finite subset of $X$) and suppose that $F$ is the family of all increasing pairs that induce $L$. According to a certain rule, we select one member of $F$ and put an arrow between the two points. 
Repeat this for all lines. Do the same for all lines, and for decreasing pairs as well. With a few exceptions, distinct arrows represent distinct lines. Whenever two arrows represent the same line, delete one of the arrows. The resulting graph $G'$ may have some isolated vertices (and indeed, the original graph $G$ may have them as well). But a careful analysis ensures that for each isolated vertex there is a nearby point with degree at least 2. Counting the degrees, we find out that there are at least the required number of arrows, and hence at least the required number of lines.

\begin{figure}
  \centering
    \includegraphics[width=0.6\textwidth]{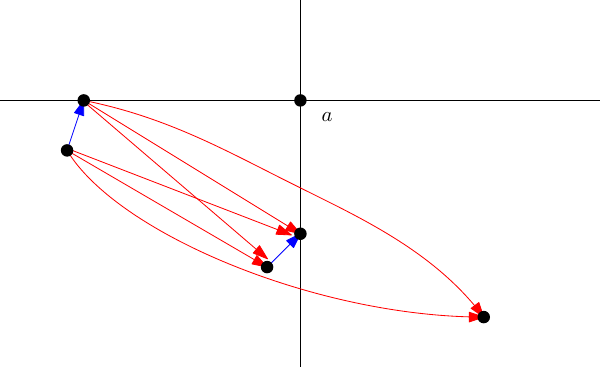}
		\caption{An example of a configuration with an isolated vertex}
		\label{fig:counterexample}
\end{figure}

 In the rest of this paper, we will always assume that $X$ contains no universal lines.

%%%%%%%%%%%%%%%%%%%%%%%%%%%%%%%%%%%%%%%%%%%%%%%%%%%%%%%%%%%5
\section{Increasing and decreasing pairs, blue and red arrows}

In this section, we introduce some necessary definitions and provide details for the construction of the graph $G$ mentioned in previous section. 

For a point $p$ in the plane, we define $x(p)$ to be the $x$-coordinate of $p$ and $y(p)$ the $y$-coordinate. 
Let $\{a,b\}$ be a pair of points in the plane. We say that it is an {\em increasing pair} if $(x(b)-x(a))\cdot (y(b)-y(a))>0$. A {\em decreasing pair} has the inequality reversed. A {\em horizontal pair} has $y(a)=y(b)$, while {\em vertical pair} satisfies $x(a)=x(b)$. 
The line defined by points $u,v$ will be denoted by $\ol{uv}$.

Let us define an equivalence $\sim_I$ on the set of increasing pairs: $\{p,q\}\sim_I \{a,b\}$ if they induce the same line. We want to select a single pair from each class of $\sim_I$. To this end we first define a partial order $\leq_I$ on $X$: we put $c \leq_I d$ if $x(c)\leq x(d)$ and $y(c)\leq y(d)$. 
Note that if $\{p,q\}$ and $\{a,b\}$ are increasing pairs and $\ol{pq} = \ol{ab}$, then the points $\{p,q,a,b\}$ form a chain in $\leq_I$. 
We will use this to define another partial order, this time on the set of increasing pairs. If $\{p,q\}$ and $\{a,b\}$ are increasing with $p\leq_I q$ and $a\leq_I b$, and $\ol{pq} =\ol{ab}$, then $\{p,q\}\leq_I^*\{a,b\}$ if either $q<_I b$, or $q=b$ and $p\geq_I a$. This is a linear order on each class of equivalence $\sim_I$. For each class of equivalence, we find the least element $\{p_0,q_0\}$ in this order. If $p_0<_I q_0$, we say that the ordered pair $(p_0,q_0)$ is a {\em blue arrow}.

For decreasing pairs we have analogous definitions. For $c,d\in X$, we put $c\leq_D d$ if $x(c)\leq x(d)$ and $y(c)\geq y(d)$. Using this, we define $\leq_I^*$ and {\em red arrows} in complete analogy to previous paragraph. Let $G$ be the directed graph on the vertex set $X$ and with edges being the red and blue arrows.

Whenever $\{a,b\}$ is an increasing pair, there is a blue arrow $(p,q)$ such that $\ol{ pq}=\ol{ ab}$. We call such ordered pair $(p,q)$ {\em the blue arrow for the line $\ol{ ab}$}. Similarly, there is a {\em red arrow for the line $\ol{ ab}$} whenever $\{a,b\}$ is a decreasing pair.

If $a$ and $b$ are points with $x(a)<x(b)$ and $y(a)<y(b)$, the {\em rectangle with corners $a$ and $b$} is the Cartesian product of the open intervals $(x(a),x(b))$ and $(y(a),y(b))$. The rectangle has {\em left boundary segment}, defined as the Cartesian product $\{x(a)\}\times (y(a),y(b))$, and similarly {\em right}, {\em top} and {\em bottom} boundary segments. 
We define the rectangle and its boundary segments analogously if $x(a)<x(b)$ and $y(a)>y(b)$. 

If  $\{a,b\}$ is a horizontal pair with $x(a)<x(b)$, a {\em strip between $a$ and $b$} is the cartesian product $(x(a),x(b))\times (-\infty,\infty)$. An {\em upper half-strip with corners $a$ and $b$} is the set $(x(a),x(b))\times (y(a),\infty)$. 
We can define {\em lower half-strip} analogously. For a vertical pair, we have a strip that is infinite in the direction of the $x$-axis, and two half-strips, {\em left} and {\em right}.

Let $a\in X$. The point $a$ defines four quadrants in the plane, numbered anticlockwise in the usual way: {\em first quadrant of $a$} consists of points $w$ such that $x(w)>x(a)$ and $y(w)>y(a)$, the {\em second quadrant of $a$} has points $w$ with $x(w)<x(a)$ and $y(w)>y(a)$, etc. Note that the quadrants in this text are considered to be open.

Finally, we call a region {\em empty} if it contains no points of $X$.

%%%%%%%%%%%%%%%%%%%%%%%%%%%%%%%%%%%%%%%%%%%%%%%%%
\section{What happens around isolated points?}\label{sec:isolated}

We will now set aside the question of whether two arrows can induce the same line and concentrate on proving that there are many arrows. In particular, if $a$ is an isolated vertex in the graph $G$, there is some nearby vertex that is the endpoint of two arrows. In this section, we prove some technical observations.

\begin{lemma}
If $a$ is not the endpoint of any arrow, then either the second or the third quadrant of $a$ is empty.
\end{lemma}

\begin{proof}
Let us suppose that this is not true. Then there is a point $b$ in the second quadrant and a point $c$ in the third quadrant of $a$. If there are more points in the second quadrant, pick one that minimizes $|y(b)-y(a)|$ and for which the horizontal line segment $(x(b),x(a))\times \{y(b)\}$ is empty. Similarly, pick $c$ in the third quadrant that minimizes $|y(c)-y(a)|$ and for which the segment $(x(c),x(a))\times \{y(c)\}$ is empty. Since $(b,a)$ is not an arrow, there is a red arrow $(p,q)$ for $\ol{ba}$, and we have $q\neq a$. Similarly, there is a red arrow $(p',q')$ for $\ol{ca}$, and $q'\neq a$. The following three paragraphs correspond to the three parts of Figure~\ref{fig:quadrantempty}.

\begin{figure}
  %\caption
  \centering
    \includegraphics[width=0.9\textwidth]{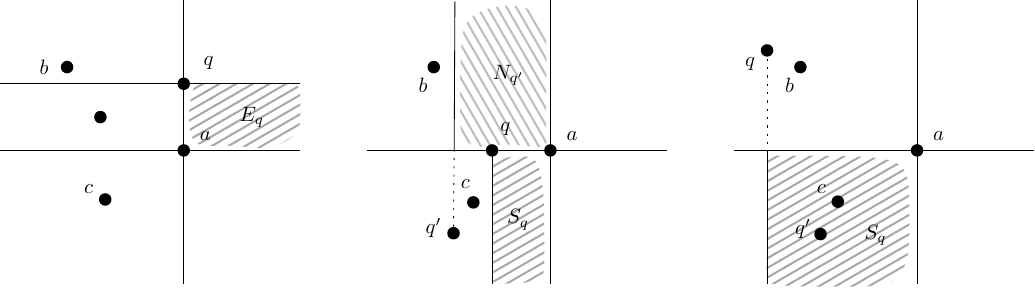}
		\caption{Second or third quadrant of an isolated point is empty.}
		\label{fig:quadrantempty}
\end{figure}

Suppose first that $x(q)=x(a)$. It follows that $y(a)<y(q)\leq y(b)$. 
Let $E_q$ be the right half-strip with corners $a,q$.  
If there is a point $z\in X$ in $E_q$, then $z\in \ol{pq} \setminus \ol{ba}$, so $E_q$ is empty. Since $\ol{aq} \neq X$, there is a point in the left half-strip with corners $a,q$. But that is a contradiction with the choice of $b$, so we may suppose that $x(q)\neq x(a)$. Similarly, $x(q')\neq x(a)$. 

Now suppose that at least one of $q,q'$ has the same $y$-coordinate as $a$. Suppose without loss of generality that $q$ has this property and if both do, then $x(q')\leq x(q)$. Let $S_q$ be the bottom half-strip with corners $a,q$, and $N_{q'}$ the upper half-strip with corners $a$ and $(x(q'),y(a))$.   
If some $z\in X$ belongs to $S_q$, then $z\in \ol{pq}\setminus\ol{ba}$, so $S_q$ is empty. Similarly $N_{q'}$ is empty. But then the upper half-strip with corners $a,q$ is empty as well (it is a subset of $N_{q'}$), so 
$\ol{qa}=X$, a contradiction. 

The only remaining case is that $q$ is in the second quadrant and $q'$ is in the third quadrant. Without loss of generality, suppose that $x(q)\leq x(q')$. Then again, the lower half-strip with corners $(x(q),y(a))$ and $a$ is empty, otherwise $\ol{pq} \neq \ol{ba}$. However, this half-strip contains $q'$, which is a contradiction. 
\end{proof}

For brevity, let us write down the following assumptions.  

\vspace{3mm}
\noindent\fbox{
	\parbox{\textwidth}{
{\em Assumptions.} Suppose that $a$ is isolated in $G$, the second quadrant of $a$ is empty, and fourth quadrant of $a$ is not empty. Let $d_1$ be a point in the fourth quadrant of $a$, picked so that the right half-strip with corners $a$ and $(x(a),y(d_1))$ is empty, and so is the line segment $(x(a),x(d_1))\times \{y(d_1)\}$. Let $d_2$ be a point in the fourth quadrant of $a$, picked so that the lower half-strip with corners $a$ and $(x(d_2),y(a))$ is empty, and so is the line segment $\{x(d_2)\}\times (y(d_2),y(a))$. 
Also, since $(a,d_1)$ is not an arrow, let $(t_1,w_1)$ be the red arrow for $\ol{ad_1}$. Similarly, let $(t_2,w_2)$ be the red arrow for $\ol{a d_2}$.  
}}
\vspace{3mm}

It may be that $d_1=d_2$. 

\begin{lemma}\label{placement}
Suppose that the Assumptions hold. We have four options for placement of $t_1$ and $w_1$:
\begin{enumerate}[(a)]
	\item $t_1$ is in the line segment $(x(a),x(d_1))\times \{y(a)\}$ and $w_1=d_1$, or
	\item $t_1$ is in the line segment $\{x(a)\}\times (y(d_1),y(a))$ and $w_1=d_1$, or
	\item $t_1$ is in the half-line $\{x(a)\}\times (y(a),\infty)$ and $w_1$ is in the line segment $(x(a),x(d_1)]\times \{y(a)\}$, or
	\item $t_1$ is in the half-line $(-\infty,x(a))\times \{y(a)\}$ and $w_1$ is in the line segment $\{x(a)\}\times [y(d_1),y(a))$.
\end{enumerate}
The four options for placement of $t_2$ and $w_2$ are analogous.
\end{lemma}

\begin{proof}
If $w_1=d_1$, then by the definition of red arrows, $a<_D t_1 <_D d_1$. Since the rectangle with corners $a$ and $d_1$ is empty, $t_1$ can only be on its top or left boundary. These are the first two options. 

Suppose that $w_1<_D d_1$. If $w_1\leq_D a$, then, since the second quadrant of $a$ is empty, $w_1$ is on one of these half-lines: $\{x(a)\}\times (y(a),\infty)$ or $(-\infty,x(a))\times \{y(a)\}$. But then, since $t_1<_D w_1$, the point $t_1$ is in the second quadrant of $a$, a contradiction. So we may suppose that $a<_D w_1 <_D d_1$, so $w_1$ is in the top or left boundary segment of the rectangle with corners $a$ and $d_1$. If it is in the top boundary, then, since $a\in \ol{t_1 w_1}$, $t_1$ is in the half-line $\{x(a)\}\times (y(a),\infty)$ and we have option (c). If it is in the left boundary, we get the last option. 
\end{proof}

%%%%%%%%%%%%%%%%%%
\begin{lemma}\label{locations}
Suppose that the Assumptions hold and that $d_1\neq d_2$. Then we have the following options for location of $(t_1,w_1)$ and $(t_2,w_2)$:

\begin{itemize}
\item $(t_1,w_1)$ follows option (a) of Lemma~\ref{placement} and $(t_2,w_2)$  option (b), or
\item $(t_1,w_1)$ follows option (a) and $(t_2,w_2)$  option (d), or
\item $(t_1,w_1)$ follows option (c) and $(t_2,w_2)$  option (b).
\end{itemize} 

Moreover, the line segments $(x(a),x(d_2)]\times \{y(a)\}$ and $\{x(a)\}\times [y(d_1),y(a))$ are both empty. 
\end{lemma}

\begin{proof}
If $t_1$ and $w_1$ are placed according to option (d), then $d_2\in \ol{t_1 w_1} \setminus \ol{a d_1}$, a contradiction. If they are placed according to option (b), then the left half-strip with corners $a$ and $t_1$ is empty: any point of $X$ that is located in it belongs to $\ol{ t_1 w_1} \setminus \ol{ a d_1}$, a contradiction. But by choice of $d_1$, the right half-strip with the same corners is empty as well, so $\ol{ a t_1}=X$. Similarly, $t_2$ and $w_2$ are not placed according to options (a) or (c). 

Out of the remaining four options for placement of the two arrows, we also rule out $(t_1,w_1)$ being placed according to (c) and $(t_2,w_2)$ according to (d): if that is the case, then $t_2\in \ol{a d_1} \setminus \ol{ t_1 w_1}$. 

Let $z$ be a point of $X$ in the line segment $\{x(a)\}\times [y(d_1),y(a))$ and suppose $(t_1,w_1)$ is placed according to (a) or (c). Then $z\in \ol{ a d_1} \setminus \ol{ t_1 w_1}$. Similarly, if $(t_2,w_2)$ is placed according to (b) or (d), the line segment $(x(a),x(d_2)]\times \{y(a)\}$ is empty.
\end{proof}

%%%%%%%%%%%%%%%%%%%%%%%%%%%%%%%%%%%%%%%%%%%%%%%%%%%5
\section{For each isolated point, there is a pair of arrows nearby}

\begin{lemma}\label{twod}
Suppose the Assumptions hold and $d_1\neq d_2$. Then $X$ contains points $z$ with $x(z)=x(a)$ and $y(z)<y(a)$. Let $r$ be such a point closest to $a$. Also, $X$ contains points $z$ with $y(z)=y(a)$ and $x(z)>x(a)$. Let $s$ be such a point closest to $a$. There is a blue arrow $(r,s)$, and a blue arrow $(r',s)$ for some point $r'$ in the fourth quadrant of $a$ such that $\{r,r'\}$ is a decreasing pair. 
\end{lemma}

\begin{proof}
By Lemma~\ref{locations} we have only three options for the locations of the arrows $(t_1,w_1)$ and $(t_2,w_2)$. In all three cases, the line segments $\{x(a)\}\times [y(d_1),y(a))$ and $(x(a),x(d_2)]\times \{y(a)\}$ are empty, and there are points in the line segments $\{x(a)\}\times [y(d_2),y(d_1))$ (either $t_2$ or $w_2$) and $(x(d_2),x(d_1)]\times \{y(a)\}$ ($t_1$ or $w_1$). 

Moreover, if $(t_2,w_2)$ is located according to option (b), then the point $t_2$ is in the open segment $\{x(a)\}\times (y(d_2),y(d_1))$. Now suppose that $(t_2,w_2)$ is located according to option (d), and $y(w_2)=y(d_2)$, that is, $\{w_2,d_2\}$ is a horizontal pair. Let $z$ be a point outside $\ol{ w_2 d_2}$. By choice of $d_2$, this point is located in the upper half-strip with corners $a$ and $(x(d_2),y(a))$. But then $z\in \ol{t_1,w_1} \setminus \ol{ a d_1}$. So if $(t_2,w_2)$ is located according to option (d), then in this case there is also a point in the open segment $\{x(a)\}\times (y(d_2),y(d_1))$. It follows that $r$ not only exists, but moreover it is located in the open segment $\{x(a)\}\times (y(d_2),y(d_1))$. Similarly, $s$ is in the open segment $(x(d_2),x(d_1))\times \{y(a)\}$. 

We will now prove that $(r,s)$ is a blue arrow. If not, then let $(p,q)$ be the blue arrow for $\ol{ r s}$. Suppose first that $q\neq s$, and hence $q<_I s$. If $q\leq_I r$, then $d_1\in \ol{ p q} \setminus \ol{ r s}$. It follows that $r<_I q <_I s$. Also, $q\neq a$ (there are no arrows into $a$), so by the choice of $d_1$, $y(q)<y(d_1)$. Also, since $q<_I s$, we have $x(q)\leq x(s)<x(d_1)$ and again, we have $d_1\in \ol {p q} \setminus \ol{ r s}$. We have shown that $q=s$. Then $r\leq_I p <_I s$, so $y(d_2)<y(r)\leq y(p)$. If $p\neq r$, then by the choice of $d_2$, $x(d_2)\leq x(p)$. It follows that $d_2\in \ol{ p q} \setminus \ol{ r s}$, a contradiction. So the only remaining option is $p=r$, and $(r,s)$ is an arrow. 

Now consider the line $\ol {d_2 s}$. If $(d_2,s)$ is not an arrow, consider the arrow $(p,q)$ for this line. The proof that $q=s$ is analogous to the one in previous paragraph, but this time $p$ is not necessarily equal to $d_2$. However, this point still has required property: $\{r,p\}$ is a decreasing pair. If not, then $r\in \ol{ p q} \setminus \ol{ d_2 s}$. Moreover, since $d_2 \leq_I p<s$, $p$ is in the fourth quadrant of $a$. Set $r'=p$.  
\end{proof}

\begin{lemma}\label{lemma:sclosest}
Suppose that Assumptions hold and $d_1=d_2$. Let $(r,s)$ be the arrow for the line $\ol{a d_1}$. According to the classification in Lemma~\ref{placement}, the points $r,s$ are located either in position (c) or (d). Moreover, there are no points on the open line segment with endpoints $a$ and $r$ and on the open line segment with endpoints $a$ and $s$.    
\end{lemma}

\begin{proof}
If $(r,s)$ is in position (a) or (b), then $\ol{ a r}=X$, a contradiction.

Now suppose that $(r,s)$ is in position (c) and let $d$ be the point $d_1=d_2$. Suppose that there is a point $u$ on the open line segment between $a$ and $s$. We have $\{r,u\}\leq_D^* \{r,s\}$, but $(r,s)$ is an arrow, so the symmetric difference of the lines $\ol{ r u}$ and $\ol{ rs}$ contains some point $z$. There are only two regions where this point might be located: the rectangle with corners $r$ and $s$, or the lower half-strip with corners $u$ and $s$. The first one is empty by the assumption $\ol{ r s}=\ol{ a d}$, and the second by choice of $d$. If there is a point $u$ on the open line segment between $a$ and $r$, then the assumption $\ol{u s}\neq \ol{r s}$ provides a contradiction similarly, this time using the assumption that the second quadrant of $a$ is empty. If $(r,s)$ is in position (d), the reasoning is completely analogous.
\end{proof}

\begin{lemma}\label{lemma:twoarrows1}
Suppose the Assumptions hold and that $d_1=d_2$. Let $(r,s)$ be the arrow for the line $\ol{ a d_1}$. Then there is another red arrow $(r',s)$ ending in $s$, where $r'\neq r$ is a point such that $\{r,r'\}$ is an increasing pair. 
\end{lemma}

\begin{proof}
The proof is very similar to the previous ones. Let $d=d_1=d_2$. Suppose $(r,s)$ is positioned according to (c) of Lemma~\ref{placement}.
Let $z$ be a point in $X$ outside $\ol{ r a}$. Since the second quadrant of $a$ is empty and $\ol{ r s}=\ol{ a d}$ and hence the rectangle with corners $r$ and $s$ is empty, $z$ can only be in the right half-strip with corners $s$ and $(x(s),y(r))$. Similarly, since $\ol{ a s}\neq X$, there is a point $c$ in the upper half-strip with corners $r$ and $(x(s),y(r))$. If $(c,s)$ is not a red arrow, let $(p,q)$ be the red arrow for $\ol{ c s}$. As before, we can argue that $q=s$ (otherwise $z\in \ol{ p q} \setminus \ol{ c s}$ and set $r'=p$. Obviously $\{r,r'\}$ is an increasing pair. In this case, $r'$ is in the first quadrant of $a$. 

If $(r,s)$ is positioned according to (d) of Lemma~\ref{placement}, then a completely analogous proof finds an $r'$ such that $\{r,r'\}$ is an increasing pair and there is an arrow $(r',s)$. This time, $s$ is below $a$ and $r'$ is in the third quadrant of $a$.
\end{proof}

According to this lemma, if the Assumptions hold and $d_1=d_2$, then either there is a point $r$ directly above $a$ (with no point between them) and $s$ directly to the right of $a$ (again, no point between them) and there are two blue arrows into $s$, one of them from $r$ and the other originating in the first quadrant of $a$. Or there is $r$ directly to the left of $a$ and $s$ directly below $a$ and again, two blue arrows into $s$, one of them from $r$ and the other originating in the third quadrant of $a$.

So far we have assumed that the fourth quadrant of $a$ is not empty. We need to see what happens if it is.

\begin{lemma}\label{lemma:twoemptyq}
Suppose that $a$ is isolated and the second and fourth quadrant are empty, while first is not. Then the half-line $(x(a),\infty)\times \{y(a)\}$ contains points of $X$. Let $s$ be such point that is closest to $a$. Also, the half-line $\{x(a)\}\times (y(a),\infty)$ contains points of $X$. Let $r$ be such  point that is closest to $a$. Then $(r,s)$ is a red arrow and there is a point $r'\neq r$ such that $(r',s)$ is a red arrow and $\{r,r'\}$ is an increasing pair. 

\end{lemma}

\begin{proof}
Let $b$ be a point in the first quadrant, chosen so that the rectangle with corners $a$ and $b$ is empty, and so are its top and right boundary segments. 

If there is only one such $b$, then the right half-strip with corners $a$ and $(x(a),y(b))$ is empty, and the upper half-strip with corners $a$ and $(x(b),y(a))$ is also empty. But since $\ol{ a b}\neq X$, there is some point $c$ in one of the half-lines $\{x(a)\}\times (y(b),\infty)$ or $(x(b),\infty)\times \{y(a)\}$. Let $(t,u)$ be the arrow for $\ol{a b}$. If $u$ is in the third quadrant of $a$, or $u=a$, then $\ol{ t u}$ contains $c$, which does not happen. So either $t$ or $u$ or both are strictly between $a,b$ in the $\leq_I$ order. It follows that there is a point $w$ with $a<_I w <_I b$ (this is either $t$ or $u$). This $w$ is either in the segment $\{x(a)\}\times (y(a),y(b)]$ or in the segment $(x(a),x(b)]\times \{y(a)\}$. Either way, $\ol{a w}=X$. 

So there are at least two points, $b_1$ and $b_2$, satisfying the condition. Pick them in such a way that the upper half-strip with corners $a$ and $(x(b_1),y(a))$ is empty and the right half-strip with corners $a$ and $(x(a),y(b_2))$ is empty.  
Again, since $(a,b_1)$ is not an arrow, there is another pair generating the same line and since this line does not contain $b_2$, one of the two points, call it $w$, is strictly between $a$ and $b_1$ in the $\leq_I$ order. In fact, it can only be located in the open line segment $\{x(a)\}\times (y(b_2),y(b_1))$.  
Let $r\in X$ be the point in this segment that is closest to $a$. 
Argue completely analogously that there is a point in the segment $(x(b_1),x(b_2))\times \{y(a)\}$ and let $s$ be such a point closest to $a$. 

The proof that $(r,s)$ is an arrow is a repetition of previous proofs. So is the proof that the arrow $(p,q)$ for the line $\ol {b_1 s}$ has $q=s$, and $\{p,r\}$ is an increasing pair. Again, set $r'=p$.
\end{proof}

If second and fourth quadrants of $a$ are empty, then one of Lemmas~\ref{twod}, \ref{lemma:twoarrows1}, or \ref{lemma:twoemptyq} applies. The resulting situations are depicted as $A,B,C$ of figure~\ref{fig:fourpics3}. 

If the third quadrant is empty and second is not, then everything above is just mirrored around the horizontal axis at $y=y(a)$. Red arrows become blue and an $r$ directly below $a$ becomes an $r$ directly above $a$. Parts $D,E,F$ of Figure~\ref{fig:fourpics3} depict the resulting situations. 

We will see in Section~\ref{sec:counting} that we may assume that there are no points $a$ such that their first, second, and fourth quadrants are all empty.

%%%%%%%%%%%%%%%%%%%%%%%%%%%%%%%%%%%%%%%%%%%%%%%%%%%%%55
\section{For each isolated point, there is a {\em unique} pair of arrows nearby}\label{sec:isolatedII}

For each isolated vertex $a$ (except those with empty first, second, and fourth quadrants), %except the one with empty first, second, fourth quadrant 
Lemma~\ref{twod}, \ref{lemma:twoarrows1}, or \ref{lemma:twoemptyq}, or their mirrored variants provide a point $s_a$ such that $\{a,s_a\}$ is a horizontal or vertical pair and there are two blue or two red arrows into $s$. Let $I$ be the set of isolated vertices $a$ such that the union of first, second, and fourth quadrant of $a$ is nonempty. Define a mapping $f:I \to X\times \{\mbox{red}, \mbox{blue}\}$ by $f(a)=(s_a,\mbox{blue})$ if the two arrows are blue, and $f(a)=(s_a,\mbox{red})$ otherwise. Note that the first component of $f$ maps $I$ into $X\setminus I$. 

\begin{lemma}\label{injective}
The mapping $f$ is injective. 
\end{lemma}

\begin{figure}
  \centering
    \includegraphics[width=0.9\textwidth]{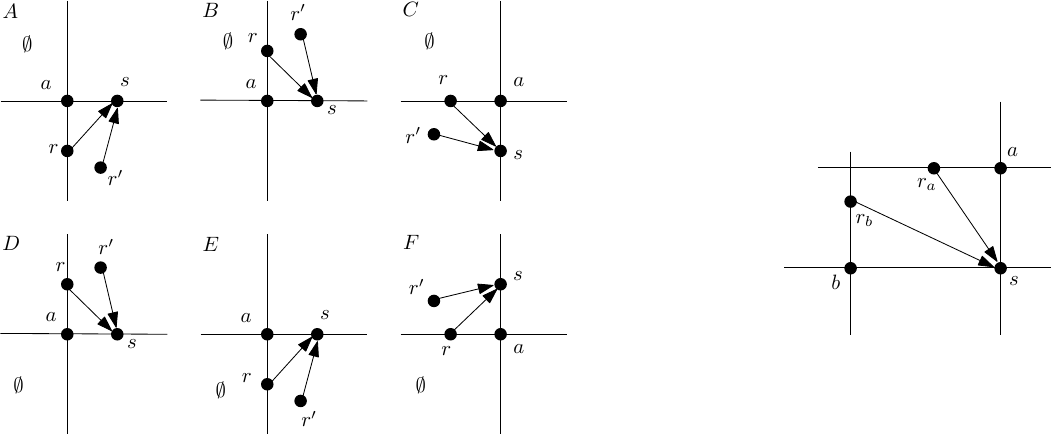}
		\caption{Four situations for isolated vertices, and an example of conflict.}
		\label{fig:fourpics3}
\end{figure}

\begin{proof}
Suppose that there is some $s$ and some $a,b\in I$ such that $f(a)=f(b)=(s,\mbox{red})$. The pair has to be increasing, with one of the red arrows (call it $(r_a,s)$) being the $(r,s)$ arrow for $a$ and the other (call it $(r_b,s)$) being the $(r,s)$ arrow for $b$. Suppose $x(b)<x(a)$. Then $a$ had to have been dealt with by the second part of Lemma~\ref{lemma:twoarrows1}. In particular, there is a $d_a$ with $x(a)<x(d_a)$ and $y(d_a)\leq y(s)$, such that the right half-strip with corners $a$ and $(x(a),y(d_a))$ is empty. (This is situation $C$ in Figure~\ref{fig:fourpics3}.) 

If the point $b$ was treated by the first part of Lemma~\ref{lemma:twoarrows1} (this is depicted as situation $B$ in the figure), the lemma provided a point $d_1=d_2$. Call it $d_b$. Since $\ol{ r_b s}=\ol{ b d_b}$, the rectangle with corners $r_b$ and $s$ is empty, and all four of its boundary segments are also empty. But second quadrant of $b$ is empty and so is the right half-strip with corners $s$ and $a$ (by the properties imposed on the point $d_a$). It follows that $\ol{ b r_b}=X$, a contradiction. 

If the point $b$ was treated by Lemma~\ref{lemma:twoemptyq} (this is also situation $B$ in the figure), then fourth quadrant of $b$ was empty. Also, the rectangle with corners $r_a$ and $s$ is empty, and so is its lower boundary segment, since $\ol{ r_a s}=\ol{ a  d_a}$.  
But then $\ol{r_a a}=X$, a contradiction.

The only remaining option is that $b$ was treated by mirrored version of Lemma~\ref{twod} (this is situation $D$). By the properties of $d_a$, the right half-strip with corners $s$ and $a$ is empty. But mirrored version of Lemma~\ref{twod} guarantees the existence of a point $d_1$ such that $y(b)<y(d_1)<y(r_b)<y(r_a)$ and $x(s)<x(d_1)$. That is, $d_1$ is in the empty region, a contradiction. 

If $f(a)=f(b)=(s,\mbox{blue})$, we reach a contradiction in an analogous way.
\end{proof}

%%%%%%%%%%%%%%%%%%%%%%%%%%%%%%%%%%%%%%%%%%%%%%%%%%%%%%%%%%%%%%%5
\section{The lines induced by arrows are (mostly) distinct}\label{sec:distinct}

\begin{lemma}\label{distinct}
	If $(p,q)$ and $(a,b)$ are two distinct arrows, then $\ol{pq}\neq \ol{ a b}$, unless the two arrows are diagonals of the same rectangle.
\end{lemma}

\begin{proof}
If the two arrows are both blue or both red, the conclusion follows from the definition. If $(p,q)$ is blue while $(a,b)$ red, then $p$ and $q$ are in the same section of the line $\ol{ab}$ (between $a$ and $b$, or one of the two tails). If they both belong to the tail of $\ol{ab}$ that contains $b$, then $\ol{pq}$ does not contain $a$, and same with the other tail. If they are both between $a$ and $b$, then the only way $\ol{pq}$ contains both $a$ and $b$ is if $a,b,p,q$ are four corners of the same rectangle. %But then $\ol{pq}=X$, a contradiction.  
\end{proof}

\begin{lemma}\label{nested}
Let $(a_1,b_1)$ and $(a_2,b_2)$ be blue arrows and $(c_1,d_1)$, $(c_2,d_2)$ red arrows such that $\ol{ a_1 b_1} = \ol{ c_1 d_1}$ and $\ol{a_2 b_2} = \ol{c_2,d_2}$. Then the two rectangles intersect in their middle parts. That is, $\{a_1,a_2\}$ and $\{b_1,b_2\}$ are decreasing pairs, while $\{c_1,c_2\}$ and $\{d_1,d_2\}$ are increasing.
\end{lemma}

\begin{proof}
Suppose that $a_1 <_I b_1 $ and $c_1 <_D d_1$. Then the third quadrant of $a_1$ is empty, including its border (except for $a_1$ itself). Same goes for the appropriate quadrants (i.e., those facing `out') with corners $b_1,c_1$ and $d_1$. The rectangle with the corners $a_2,d_2,b_2,c_2$ therefore fits between the horizontal lines $y=y(a_1)$ and $y=y(c_1)$, or between the vertical lines $x=x(a_1)$ and $x=x(d_1)$. The same is true if the roles of the two rectangles are exchanged, so the only placement possible is the one described.
\end{proof}

\begin{figure}
  \centering
    \includegraphics[width=0.9\textwidth]{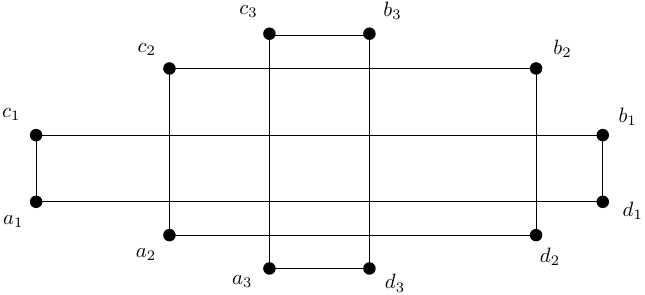}
		\caption{Blue and red arrows generating the same lines.}
		\label{fig:nested}
\end{figure}

Let $(a_i,b_i)$ with $a_i<_I b_i$ for $i=1,\dots,k$ be blue arrows and $(c_i,d_i)$ with $c_i<_D d_i$ be red arrows such that $\ol{ a_i b_i}=\ol{ c_i d_i}$ for all $i$. It follows from Lemma~\ref{nested} that we can assume that the points are numbered in such a way that $a_1<_D a_2<_D \dots <_D a_k$, $b_1>_D \dots >_D b_k$, $c_1<_I \dots <_I c_k$ and $d_1>_I \dots >_I d_k$. See Figure~\ref{fig:nested}. 

If two arrows $(a_i,b_i)$ and $(c_i,d_i)$ induce the same line, we delete one of them. We use the following lemma to show that we still have enough arrows.

\begin{lemma}\label{counting}
Let $a_i$, $b_i$, $c_i$ and $d_i$ be points as described above. For each $i$ there is an arrow starting in $a_i$, other than $(a_i,b_i)$. Also for each $c_i$, there is an arrow starting in it, other than $(c_i,d_i)$.
% of i<k? Check.
\end{lemma}

\begin{proof}
Let $u_1,\dots,u_{\ell}$ be all points in $X$ whose (open) third quadrant is empty. Every two form either a horizontal, a vertical, or a decreasing pair. We can assume they are numbered so that $u_1<_D u_2<_D \dots <_D u_{\ell}$. Since the sequence contains all points with the property, the only points of $X$ in the closed rectangle (or a segment, if the pair is horizontal or vertical) $[x(u_j),x(u_{j+1})]\times [y(u_{j+1}),y(u_j)]$ are $u_j$ and $u_{j+1}$. 

The third quadrants of points $a_1,\dots,a_k$ are empty (even including their borders), so they fall somewhere in the list of the $u_i$'s. Note that if the third quadrant of $c_1$ is empty, then $c_1$ belongs among the $u_t$'s, so $a_1\neq u_1$. If the third quadrant of $c_1$ is not empty, then again we find some other $u_t$ preceding $a_1$ and again $a_1\neq u_1$. Similarly, $a_k\neq u_{\ell}$. The point is that for each $i$, we find $j$ such that $a_i=u_j$, and for each such $j$, both $u_{j-1}$ and $u_{j+1}$ exist.

Fix $i$ with $1< i\leq k$ and let $j$ be such that $u_j=a_i$. The boundary of the third quadrant of $u_j$ is empty, so $\{u_{j-1},u_j\}$ is either a decreasing or a vertical pair. Similarly, $\{u_j,u_{j+1}\}$ is either decreasing or horizontal.

\begin{itemize}
\item Suppose first that both pairs are decreasing. Pick an $s$ with $1\leq s\leq k$ such that $s\neq i$. There are some points in the rectangle $(x(u_j),x(b_s)]\times (y(u_j),y(b_s)]$ (for instance, $b_s$ itself). From these, let us choose a $q_j$ such that the rectangle $(x(u_j),x(q_j)]\times (y(u_j),y(q_j)]$ only contains $q_j$. Note that the location of $q_j$ guarantees that $q_j \neq b_i$. 

If $(u_j, q_j)$ is not an arrow, let $(p,q)$ be the blue arrow assigned to the line $\ol{ u_j q_j}$. The third quadrant of $u_j$ is empty, including its boundary, so we have $u_j \leq_I p<_I q \leq_I q_j$. The only possible placement is that $q=q_j$ and $p$ is on the horizontal segment $(x(u_{j+1}),x(q_j))\times \{y(u_j)\}$, or on the vertical segment $\{x(u_j)\}\times (y(u_{j-1}),y(q_j))$. 
% **** definovat arrow of the line
In both cases, either $u_{j-1}$ or $u_{j+1}$ belong to $\ol{ p q} \setminus \ol{ u_j q_j}$, which is a contradiction.

\item If $\{u_{j-1},u_j\}$ is decreasing and $\{u_j,u_{j+1}\}$ horizontal, then we may pick $q_j$ such that the rectangle $(x(u_j),x(q_j)]\times (y(u_j),y(q_j)]$ only contains $q_j$, with the additional property that $q_j$ belongs to the vertical strip $(x(u_j),x(u_{j+1}))\times (y(u_j),\infty)$. There are such points, since $\ol{ u_j u_{j+1}} \neq X$ and the third quadrant of $u_{j+1}$ is empty. 

We have $x(u_j)<x(q_j)<x(u_{j+1})\leq x(d_k) =x(b_k)\leq x(b_i)$. 
This shows that $q_j\neq b_i$.

If $(p,q)$ is, as before, the blue arrow for $\ol{u_j q_j}$, then $q=q_j$ and the only placement for $p$ is on the segment $\{x(u_j)\}\times (y(u_{j-1}),y(q_j))$. But then $u_{j-1}\in \ol{ p q} \setminus \ol{ u_j q_j}$. 

\item If $\{u_{j-1},u_j\}$ is vertical and $\{u_j,u_{j+1}\}$ decreasing, then by similar reasoning pick a $q_j$ in the horizontal strip $(x(u_j),\infty)\times (y(u_j),y(u_{j-1}))$ with the property that $(x(u_j),x(q_j)]\times (y(u_j),y(q_j)]$ only contains $q_j$. As before, $q_j\neq b_i$. Again, let $(p,q)$ be the blue arrow assigned to the line $\ol{ u_j q_j}$ and observe that $q=q_j$ and the only placement for $p$ is on the segment $(x(u_{j+1}),x(q_j))\times \{y(u_j)\}$. But then $u_{j+1}\in \ol{ p q} \setminus \ol{ u_j q_j}$. 

\item Now suppose that $\{u_{j-1},u_j\}$ is vertical and $\{u_j,u_{j+1}\}$ horizontal. Again, pick $s\neq i$. Let $(p,q)$ be the blue arrow for the line $\ol{ u_j b_s}$. If $p=u_j$, then we have found an arrow from $u_j=a_i$ somewhere else than $b_i$ and we are done. 

If $p\neq u_j$, then $p$ is in the region $A:=[x(u_{j+1},x(b_s)]\times [y(u_{j-1}),y(b_s)]$ (otherwise $u_{j-1}$ or $u_{j+1}$ belong to $\ol{ p q}\setminus \ol{ u_j q_j}$). Let $p'$ be a point in $A$ such that $[x(u_{j+1}),x(p')]\times [y(u_{j+1}),y(p')]$ only contains $p'$. We claim that $(u_j,p')$ is an arrow. If not, let $(a,b)$ be the blue arrow for $\ol{ u_j p'}$. Since the third quadrant of $u_j$ is empty and so is its boundary, we have $u_j\leq_I a<_I b\leq_I p'$. Since $u_{j-1}$ and $u_{j+1}$ are in $\ol{ a b}$, $b\in A$ and hence $b=p'$. Similarly, $a=u_j$. So $(u_j,p')$ is an arrow. 
\end{itemize}

We can use completely analogous arguments for the sequence $c_1,\dots, c_k$ and find an arrow from each of them to a point different than the corresponding $d_i$.
\end{proof}

%%%%%%%%%%%%%%%%%%%%%%%%%%%%%%%%%%%%%%%%%%%%%%%%%%%%%%%
\section{Counting degrees---there are many lines}\label{sec:counting}

\begin{proof}[Proof of Theorem~\ref{main}.]
First, suppose that $X$ has more than one vertex whose first, second, and fourth quadrant are empty. Pick any two such vertices, say $a$ and $b$. They are a horizontal or vertical pair, and the line $\ol{ a b}$ contains all of $X$. So we may suppose that there is at most one such point, call it $a$. If there are any points on the half-line $\{x(a)\}\times (y(a),\infty)$ or on the half-line $(x(a),\infty)\times \{y(a)\}$, then again any of these points together with $a$ defines a line that contains all of $X$. So all points of $X$ belong to the third quadrant of $a$ and its boundary. If there is a point $c\in X$ for which second, third, and fourth quadrant is empty, then by similar reasoning, all of $X$ belongs to the first quadrant of $c$ and its boundary. But then $\ol{ a c}=X$, a contradiction. So the existence of $a$ prevents the existence of $c$. If there is such $a$, we may create a set $X'$ by reflecting $X$ in the origin. The resulting $X'$ has equal number of lines as $X$, but it has no points with first, second, and fourth quadrant empty. The reasoning in the following paragraphs would then be applied to this $X'$. So we may without loss of generality assume that such $a$ does not exist.

Let us consider the directed graph $G$ on the vertex set $X$, where the edges are the red and blue arrows. This graph might have some isolated vertices. In Section~\ref{sec:isolated}, we defined a mapping $f$ that assigns a pair $(s_a,\mbox{red})$ or $(s_a,\mbox{blue})$ to each isolated $a$, where $s_a$ is a point located close to $a$ and the second entry denotes the color of the two arrows ending in $s_a$ that are assigned to this $a$. We have also shown that this mapping is injective. 

Let $u$ be the number of isolated vertices. We have $u$ pairs $(a,s_a)$, and each pair has at least two arrows ending in $s_a$ associated with it. The pairs may overlap (namely we may have $s_a=s_b$ for distinct isolated vertices $a$ and $b$), but we have already shown that the pairs of arrows are pairwise disjoint. If $I$ is the set of isolated vertices, the set $U:=I\cup \{s_a; a\in I\}$ has at most $2u$ vertices, and 
\[ \sum_{v\in U} \deg^+(v)\geq 2u.
\] 
As before, let $(a_i,b_i)$ and $(c_i,d_i)$ for $i\in \{1,\dots,k\}$ be the coinciding red and blue arrows. Let $V:=\{a_1,\dots,a_k,c_1,\dots,c_k\}$. We have shown in Section~\ref{sec:distinct} that each vertex in $V$ has outdegree at least 2. We have $|U\cup V|\leq 2u+2k$, so $W:=X\setminus (U\cup V)$ has size at least $n-2u-2k$. The vertices outside $U$ and $V$ have at least one arrow starting or ending in them. We have 
\begin{eqnarray*}
2|E(G)|&=&\sum_{v\in X} \deg^+(v) + \sum_{v\in X} \deg^-(v) \\
&\geq& \sum_{v\in U} \deg^+(v) + \sum_{v\in V} \deg^-(v) + \sum_{v\in W} \deg^+(v) + \sum_{v\in W} \deg^-(v) \\
&\geq& 2u+2\cdot 2k + (n-2u-2k) \\
&=& n+2k.
\end{eqnarray*}
Now delete the arrow $(a_i,b_i)$ for each $i=1,\dots,k$. This modified graph $G'$ has $2|E(G')|\geq n$. It follows that $G'$ has at least $n/2$ arrows, each corresponding to a distinct line of the original metric space. 
\end{proof}

%%%%%%%%%%%%%%%%%%%%%%%%%%%%%%%%%%%%%%%%%%%%%%%%%%%%%%%%%%%%%%%5
\section{$L_{\infty}$ metric}

Let us briefly consider the $L_{\infty}$ metric on $\mathbb{R}^d$, defined by 
\[ d((u_1,\dots,u_d),(v_1,\dots,v_d))= \max_{1\leq j \leq d} |u_j-v_j|.
\]

Whenever we have a finite set in the plane with the $L_{\infty}$ metric, rotating the plane by 45 degrees transforms the $L_{\infty}$-lines into $L_1$-lines. The following theorem is an easy consequence of Theorem~\ref{main}.

\begin{theorem}\label{main2}
Let $X$ be a set of $n$ points in the plane with the $L_{\infty}$ metric. If there is no universal line, then $X$ induces at least $\lceil n/2\rceil$ lines.
\end{theorem}

This may be of interest because any finite metric space can be regarded as an $L_{\infty}$ space in $\mathbb{R}^d$ for some $d$. However, the required $d$ may be very high, which makes handling it difficult.

\vspace{3mm}

{\bf Acknowledgment.} I want to thank Va\v{s}ek Chv\'{a}tal for his many insights into the topic of lines in metric spaces, as well as useful remarks concerning this paper.
%%%%%%%%%%%%%%%%%%%%%%%%%%%%%%%%%%%%%%5

%\bibliographystyle{siam}
%\bibliography{ida}

\end{document}